\documentclass[12pt]{article} \usepackage{amstex}

\usepackage{amssymb} \usepackage{amscd} 
\usepackage{theorem} \usepackage{a4} 

\def\Talpha#1{\vbox{\ialign{##\crcr
$\alpha$\crcr\noalign{\kern2pt\nointerlineskip}
$\hfil\displaystyle{#1}\hfil$\crcr}}}

\def\Onabla#1{\vbox{\ialign{##\crcr
$\,\scriptstyle{0}$\crcr\noalign{\kern2pt\nointerlineskip}
$\hfil\displaystyle{#1}\hfil$\crcr}}}

\def\cala{{\cal A}}

\def\calc{{\cal C}}

\def\bbbone{\mbox{\rm 1\hspace {-.6em} l}}

\def\gder{{\mbox{Der}}}

\newtheorem{lemma}{LEMMA} \newtheorem{proposition}{PROPOSITION}

\begin{document}

\baselineskip=0.7cm \begin{center} {\Large\bf SHADOW OF NONCOMMUTATIVITY}
\end{center} \vspace{0.75cm}

\begin{center} Michel DUBOIS-VIOLETTE \\
John MADORE\\
\vspace{0.3cm} {\small Laboratoire de Physique Th\'eorique et Hautes
Energies\footnote{Laboratoire associ\'e au Centre National de la Recherche
Scientifique - URA D0063}\\ Universit\'e Paris XI, B\^atiment 211\\ 91 405
Orsay
Cedex, France\\ flad$@$qcd.th.u-psud.fr, madore$@$qcd.th.u-psud.fr}\\
and\\
Richard KERNER\\
\vspace{0.3cm} {\small Laboratoire de Physique Th\'eorique\\
Gravitation et Cosmologie Relativistes\footnote{Laboratoire associ\'e au
Centre
National de la Recherche Scientifique - URA D0769}\\
Tour 22/12 - 4\`eme \'etage - Bo\^\i te 142\\
4, place Jussieu - 75005  Paris, France\\
rk$@$ccr.jussieu.fr}

\end{center} \vspace{1cm}

\begin{center} \today \end{center}

\vspace {1cm}

\noindent L.P.T.H.E.-ORSAY 96/06\\

\newpage
\begin{abstract}
We analyse the structure of the $\kappa=0$ limit of a family of algebras
${\cal A}_\kappa$ describing noncommutative versions of space-time, with $\kappa$
a parameter of noncommutativity. Assuming the Poincar\'e covariance of the
$\kappa=0$ limit, we show that, besides the algebra of functions on Minkowski
space, ${\cal A}_0$ must contain a nontrivial extra factor ${\cal A}^I_0$ which is
Lorentz covariant and which does not commute with the functions whenever it is
not commutative. We give a general description of the possibilities and analyse
some representative examples.
\end{abstract}

\section{Introduction} In this paper we wish to give a heuristic analysis of
the
structure of a noncommutative space-time, as might arise from the
interplay of quantum theory with gravity, in the limit where the space-time is
the usual commutative Minkowski space, that is in
the limit when the Planck length tends to zero (when gravity is switched off).
We shall show that in this limit, one necessarily captures a non-trivial extra
factor in the limit algebra
besides the algebra $\calc(M)$ of functions on Minkowski space $M$. In view of
what is known of particle
interactions, it is natural to expect that this extra factor has something to
do
with gauge theory. We have been aware of this conclusion for some time; it was
indeed one of our main motivations in studying gauge theory over the
noncommutative
algebra of matrix-valued functions \cite{dvkm}; the model for the extra factor
being
then described by the matrix algebra. We shall see however that a matrix
algebra
is too small to be the extra factor and that furthermore a noncommutative extra
factor cannot commute with the algebra $\calc(M)$.\\

We recall briefly the arguments which suggest that the interplay between
quantum
theory and gravitation leads to a noncommutative space-time. There is first an
old semi-classical argument which is recalled in \cite{DFR}, showing that
localization looses its meaning at distances of the order of the Planck length
$\lambda_p$. The argument is that, because of quantum theory, in order to
localize an event to within $\Delta x^\mu\sim a$, one needs to transfer an
amount
of energy of order $1/a$ and that then, in view of general relativity, if $a$
is
too small, say $a<\lambda_p$, the energy would create a ``black hole". In fact
this
semi-classical argument can be made more precise \cite{DFR} and leads to
limitations
of the form $\Delta x^0\sum \Delta x^k\gtrsim \lambda^2_p$ and
$\displaystyle{\sum_{k<\ell}} \Delta x^k\Delta x^\ell \gtrsim \lambda^2_p$ in
order to avoid the above phenomena. This suggests that the $x^\mu$ do not
commute, that the space-time is noncommutative. It is worth noticing here
that this is not the only conclusion. For instance, one can argue that since in
this argument the $x^\mu$ have a length scale, then the metric must enter
somewhere and that
the above uncertainty relations can be consequences of the quantization of the
metric instead of a ``quantization" of space-time itself.\\
There is however a
second argument of a different nature. In  classical general
relativity, one solves locally Einstein's equations and one extends maximally
the
solutions. Doing this one obtains possibly a space-time carrying a non-trivial
topology. This extension is not just a mathematical artifact but is physically
relevant. Consider, for instance, the analogue of the classical self-energies
of
point charges. More generally, this has to be taken into account in order to
use
the old physical idea \cite{DS} that gravitation acts as an ultraviolet
``regularization" because of the fact that it is attractive and has also its
own
energy-momentum density source. A problem arises when one quantizes the
gravitational field. The topology of space-time must then be sensitive to the
states of
the field. This suggests again that the functions on space-time should be
replaced by
elements of a noncommutative $\ast$-algebra acting on the same states as the
quantized gravitation field.\\

In short there are several arguments suggesting that the space-time
becomes noncommutative on the scale of Planck length. Although this is
certainly not the only
possibility, it is worth studying the consequences which follow from such an
assumption when the gravitational interaction is switched off. It is the aim of
next section to analyze the ``shadow" of such a noncommutativity in  the limit
when one recovers the Poincar\'e-covariant physical theory in the usual
Minkowski
space. This analysis is sharpened in Section 3 where it is pointed out that the
limit has a structure of crossed product and  the Hochschild 2-cocycle
associated with the deformation is identified as a group cocycle of the dual of
the group of space-time translations. Then, in Section 4, we shall discuss in
this context various recent
proposals and, in particular, the generalizations of gauge theory using
noncommutative differential calculi \cite{dvkm}, \cite{C}, \cite{coq}.

\section{The commutative space-time limit}

Let ${\cal A}_\kappa$ be a one-parameter family of associative algebras. Here
we
think of $\kappa$ as being the gravitational constant $\lambda_p^2$ and we
think
of  $\cala_\kappa$ as being a noncommutative version of the functions on
space-time. Technically, we use  the framework of formal deformation theory of
associative algebras \cite{G}, \cite{BFLS}. This means that we assume that, as
vector
spaces, all the $\cala_\kappa$ coincide with a fixed vector space $E$ and that,
for $f,g\in E$, one can expand their product $(fg)_\kappa$ in $\cala_\kappa$ as
\begin{equation} (fg)_\kappa=fg+\kappa c(f,g)+o(\kappa^2) \label{(eq1)}
\end{equation} where $fg=(fg)_0$ is the product in $\cala_0$. We also assume
that
there is a distinguished element, $\bbbone \in E$, which is the unit for each
algebra $\cala_\kappa$. The bilinear map $(f,g)\mapsto c(f,g)$ is then a
normalized Hochschild 2-cocycle of $\cala_0$ with values in $\cala_0$ \cite{G}.
In terms of commutators, (\ref{(eq1)}) yields \begin{equation}
[f,g]_\kappa=[f,g]-i\kappa\{f,g\}+o(\kappa^2) \label{(eq2)} \end{equation}
where
the {\sl bracket} $(f,g)\mapsto \{f,g\}$ is defined by \begin{equation}
\{f,g\}=i(c(f,g)-c(g,f)) \label{eq3} \end{equation} and where $[f,g]_\kappa$
and
$[f,g]$ denote the commutator in $\cala_\kappa$ and in $\cala_0$ respectively.
The first order
in $\kappa$ of the identity $[h,(fg)_\kappa]_\kappa=([h,f]_\kappa
g)_\kappa+(f[h,g]_\kappa)_\kappa$ yields the condition \begin{equation}
i([h,c(f,g)]-c([h,f],g)-c(f,[h,g]))=f\{h,g\}-\{h,fg\}+\{h,f\}g \label{eq4}
\end{equation} This implies that, if $h$ is an element of the center
$Z(\cala_0)$ of $\cala_0$, then the endomorphism $\delta_h$ of $\cala_0$
defined
by $\delta_h(f)=\{h,f\}$ is a derivation of $\cala_0$. The center $Z(\cala_0)$
of
$\cala_0$ is stable under the derivations of $\cala_0$ and therefore
$Z(\cala_0)$ is stable under the bracket (\ref{eq3}). On the other hand, if $f$, $g$ and $h$ are elements of $Z(\cala_0)$,
then the lowest non-trivial order in $\kappa$ (i.e. the second order) of
$[[f,g]_\kappa,h]_\kappa=[[f,h]_\kappa,g]_\kappa+[f,[g,h]_\kappa]_\kappa$
yields the condition
\begin{equation}
\{\{f,g\},h\}=\{\{f,h\},g\}+\{f,\{g,h\}\} 
 \label{eq5} \end{equation}
So one can summarize the above
discussion by the following proposition:

\begin{proposition} The center $Z(\cala_0)$ of $\cala_0$ is a (commutative)
Poisson algebra with Poisson bracket given by (\ref{eq3}) and one defines a linear mapping $z\mapsto \delta_z$ of $Z(\cala_0)$ into the Lie algebra
$\gder(\cala_0)$ of all derivations of $\cala_0$ by setting
$\delta_z(f)=\{z,f\}$
for $z\in Z(\cala_0)$ and $f\in \cala_0$ which satisfies $\delta_{zz'}=z\delta_{z'}+z'\delta_z$ for $z,z'\in Z(\cala_0)$. \end{proposition}
This result is known, it is the first part of Proposition 1.2 of \cite{RVW}
(see also \cite{DKP}).

We wish to represent the noncommutative analogue of real functions by hermitian
elements. This leads us to add the following reality condition to the above
general structure. We assume that the $\cala_{\kappa}$ are complex
$\ast$-algebras such that all the underlying complex involutive vector spaces
coincide. Thus the vector space $E$ must be a complex vector space equipped
with an antilinear involution, $f\mapsto f^\ast$, and a distinguished hermitian
element $\bbbone=\bbbone^\ast$. The parameter $\kappa$ is real (consistently)
and one has $(fg)_{\kappa}^\ast=(g^\ast f^\ast)_{\kappa}$ and
$(f\bbbone)_{\kappa}=(\bbbone f)_{\kappa}=f$. It follows that the normalized
cocycle $c$ satisfies $c(f,g)^\ast=c(g^\ast,f^\ast)$,  which implies that the
bracket (3) is real, i.e. $\{f,g\}^\ast=\{f^\ast,g^\ast\}$. Therefore, {\sl the
set $Z_{\Bbb R}(\cala_0)$ of all hermitian elements of the center $Z(\cala_0)$
of $\cala_0$ is a real (commutative) Poisson algebra and the map $z\mapsto
\delta_z$ induces a  real linear mapping of $Z_{\Bbb
R}(\cala_0)$ into the real Lie algebra $\gder_{\Bbb R}(\cala_0)$ of
all hermitian derivations of $\cala_0$.}\\

We now return to our specific problem. The $\cala_{\kappa}$ are noncommutative
versions of the algebra of functions on space-time and we wish to recover
Poincar\'e-invariant physics on Minkowski space $M$ in the limit $\kappa=0$. We
thus assume that the Poincar\'e group $\mathfrak P$ acts by
$\ast$-automorphisms $(\Lambda, a)\rightarrow \alpha_{(\Lambda,a)}$ on
$\cala_0$, that $\cala_0$ contains as a $\ast$-subalgebra the (commutative)
algebra $\calc(M)$ of (smooth) functions on Minkowski space and that the action
of $\mathfrak P$ on $\cala_0$ induces its usual action on $\calc(M)$ associated
with the corresponding transformations of $M$. We now argue that $\calc(M)$
cannot be the whole $\cala_0$. In fact, assume that $\calc(M)$ is equal to
$\cala_0$. Then, in view of Proposition 1, there is a Poisson bracket on $M$.
This Poisson bracket is non trivial since we have assumed that the
$\cala_{\kappa}$ are noncommutative. On the other hand there does not exist a
non trivial Poincar\'e invariant Poisson bracket on $M$. It seems unreasonable
to us that the Poincar\'e invariance be broken at the first order in $\kappa$.
In fact, at this order, we expect a Poincar\'e invariant theory involving a
spin-2 field linearly coupled to the other fields. Once this hypothesis is
accepted, it follows that the inclusion $\calc(M)\subset \cala_0$ must be a
strict one;  the $\kappa=0$ limit $\cala_0$ of the $\cala_{\kappa}$ must
contain an extra factor. It follows also that the normalized 2-cocycle $c$ of
$\cala_0$ defined by (1) is Poincar\'e invariant, so one has the condition

\begin{equation}
\alpha_{(\Lambda,a)}(c(f,g))=c(\alpha_{(\Lambda,a)}(f),\alpha_{(\Lambda,a)}(g)),\ \forall (\Lambda,a)\in \mathfrak P,\ \forall f,g\in \cala_0
\label{eq6}
\end{equation}
which implies the invariance of the bracket (3).\\

Let $x^\mu\in \calc(M)$ be minkowskian coordinates. Then the algebra $\calc(M)$
is generated by the $x^\nu$ and the action of the Poincar\'e group on it is
given by
\begin{equation}
\alpha_{(\Lambda,a)}x^\mu=\Lambda^{-1\mu}_{\phantom{-1}\nu}
(x^\nu-a^\nu\bbbone)
\label{eq7}
\end{equation}
By choosing an origin, one can identify $\calc(M)$ with the Hopf algebra of
functions on the group of translations. Since $\calc(M)$ is a subalgebra of
$\cala_0$, the latter is (in particular) a bimodule over $\calc(M)$.
Furthermore, by restricting attention to the action of translations, $\cala_0$
is in fact a bicovariant bimodule over the algebra of functions on the group of
translations \cite{W}, \cite{B}. By standard arguments \cite{W}, \cite{B},
$\cala_0$ is isomorphic as left $\calc(M)$-module to $\calc(M)\otimes
\cala^I_0$ where $\cala^I_0$ denotes the subalgebra of translationally
invariant elements of $\cala_0$: $\cala^I_0=\{f\in\cala_0\vert
\alpha_{(1,a)}f=f,\ \forall a\}$. In fact $\cala_0$ is isomorphic to
$\calc(M)\otimes \cala^I_0$ as $(\calc(M),\cala^I_0)$-bimodule. Thus in order
to recover the complete structure of algebra of $\cala_0$ in the representation
$\calc(M)\otimes \cala^I_0$, it is sufficient to describe the right
multiplication by elements of $\calc(M)$ of elements of $\cala^I_0$. For each
minkowskian coordinate $x^\mu$, it follows from (\ref{eq7}) that $\cala_0^I$ is
stable under the derivation $f\mapsto ad(x^\mu)(f)=[x^\mu,f]$ and therefore one
has $fx^\mu=x^\mu f-ad(x^\mu)(f),\ \forall f\in \cala^I_0$ which can be written
in the representation $\cala_0=\calc(M)\otimes \cala^I_0$
\begin{equation}
fx^\mu=x^\mu\otimes f-\bbbone \otimes ad(x^\mu)(f),\ \forall f\in \cala^I_0
\label{eq8}
\end{equation}
{}From this one deduces the right multiplication by elements of $\calc(M)$ i.e.
the  right $\calc(M)$-module structure of $\calc(M)\otimes \cala^I_0$, which is
in fact isomorphic to $\cala^I_0\otimes \calc(M)$ as right $\calc(M)$-module.
In the following, we shall denote by $X^\mu$ the four commuting derivations of
$\cala^I_0$ induced by the $ad(x^\mu)$. The algebra $\cala_0^I$ is invariant
under the automorphisms $\alpha_{(\Lambda,0)}$. We shall denote by $\alpha^I:
\Lambda \mapsto \alpha^I_\Lambda$, the corresponding homomorphism of the
Lorentz group into the group $\mbox{Aut}(\cala^I_0)$ of all
$\ast$-automorphisms of $\cala^I_0$. \\

One can thus summarize the above discussion by the following presentation of
$\cala_0$. One starts with a unital $\ast$-algebra $\cala^I_0$ equipped with
four commuting antihermitian derivations $X^\mu$ and an action $\Lambda\mapsto
\alpha^I_\Lambda$ of the Lorentz group by automorphisms of $\cala^0_I$ such
that
\begin{equation}
\alpha^I_\Lambda\circ X^\mu=\Lambda^{-1\mu}_{\phantom{-1}\nu} X^\nu\circ
\alpha^I_\Lambda
\label{eq9}
\end{equation}
and $\cala_0$ is ``generated" as unital $\ast$-algebra by $\cala^I_0$ and four
hermitian elements $x^\mu$ with relations $x^\mu x^\nu=x^\nu x^\mu$ and $x^\mu
f=f x^\mu+X^\mu(f),\  f\in\cala^I_0$. We put quotes on the word generated in
the above sentence because we do not pay attention here to functional analysis
aspects, for instance appropriate completions, etc.
 The Poincar\'e group acts on $\cala_0$ by the action $\alpha_{(\Lambda,a)}$ on
$\calc(M)$ defined by (\ref{eq7}) and by
$\alpha_{(\Lambda,a)}(f)=\alpha^I_\Lambda(f)$ for $f\in\cala^I_0$.\\

We have assumed that the bracket (\ref{eq3}) is not identically zero on
$\calc(M)$. This implies that the $c^{\mu\nu}=c(x^\mu,x^\nu)$ do not all
vanish. On the other hand, it follows from (\ref{eq7}) that $c^{\mu\nu}$ are
elements of $\cala^I_0$ and (\ref{eq6}) implies then that one has
\begin{equation}
\alpha^I_\Lambda(c^{\mu\nu})=\Lambda^{-1\mu}_{\phantom{-1}\alpha}
\Lambda^{-1\nu}_{\phantom{-1}\beta} c^{\alpha\beta}
\label{eq10}
\end{equation}
Thus the homomorphism $\alpha^I:{\mathfrak L}\rightarrow \mbox{Aut}(\cala^I_0)$
of the Lorentz group ${\mathfrak L}$ into the group of $\ast$-automorphisms of
$\cala^I_0$ {\sl is never trivial}. This places severe restrictions on the
structure of the extra factor $\cala^I_0$. For instance, $\cala^I_0$ cannot be
a finite-dimensional matrix algebra because on such an algebra all
automorphisms are inner and, on the other hand, it is known that the Lorentz
group has no nontrivial finite-dimensional unitary representation. The same
argument shows that, if $\cala^I_0$ admits an injective (eventually unbounded)
$\ast$-representation in a Hilbert space, then $\cala^I_0$ must be
infinite-dimensional.

\section{Crossed-product structure of $\cala_0$}

There is another useful presentation of $\cala_0$ which is a sort of
exponentiation of the previous one and which can be used for a
functional-analytic development of the framework. Let $\mathfrak  T^\ast$ be
the space $\Bbb R^4$ considered as the dual of the group ${\mathfrak T}$ of
translations of the Minkowski space $M$. Instead of taking the $x^\mu$ as
generators of $\calc(M)$ we now choose the exponentials $u(k)=\exp (ik_\mu
x^\mu)$ for $k\in {\mathfrak T}^\ast$. One has
\begin{equation}
u(0)=\bbbone,\ \  u(k)u(\ell)=u(k+\ell),\ \  u(k)^\ast=u(-k)=u(k)^{-1}
\label{eq11}
\end{equation}
and (\ref{eq7}) now yields with $(\Lambda k)_\mu=k_\nu
\Lambda^{-1\nu}_{\phantom{-1}\mu}$
\begin{equation}
\alpha_{(\Lambda,a)} u(k) = u(\Lambda k) \exp (-i(\Lambda k)_\mu a^\mu)
\label{eq12}
\end{equation}
If $f\in \cala^I_0$, then $\tau_k(f)=u(k)f\ u(-k)$ is also in $\cala^I_0$ and
$k\mapsto \tau_k$ is a homomorphism of the additive group ${\mathfrak T}^\ast$
into the group of $\ast$-automorphisms of $\cala^I_0$. This homomorphism is in
fact the exponential version of the derivations $X^\mu$, ($\tau_k=\exp(i k_\mu
X^\mu))$, and (\ref{eq9}) is replaced by
\begin{equation}
\alpha^I_\Lambda \circ \tau_k=\tau_{\Lambda k}\circ \alpha^I_\Lambda
\label{eq13}
\end{equation}
The $\ast$-algebra $\cala_0$ is generated by the $\ast$-algebra $\cala^I_0$ and
the $u(k)$, $k\in{\mathfrak T}^\ast$ with the  relations (\ref{eq11}) and the
relation
\begin{equation}
u(k)f=\tau_k(f) u(k),\ \forall f\in \cala^I_0\ \mbox{and}\ \forall k\in
{\mathfrak T}^\ast
\label{eq14}
\end{equation}
The action of the Poincar\'e group on $\cala_0$ is given by (\ref{eq12}) and by
\begin{equation}
\alpha_{(\Lambda,a)}(f)=\alpha^I_\Lambda(f),\ \forall f\in \cala^I_0
\label{eq15}
\end{equation}
The consistency is ensured by the relation (\ref{eq13}). We can summarize the
above discussion by the following proposition:
\begin{proposition}
The $\ast$-subalgebra $\cala^I_0$ of translationally invariant elements of
$\cala_0$ is equipped with an action $\tau$ of the dual group $\mathfrak
T^\ast$ of the group $\mathfrak T$ of translations by $\ast$-automorphisms and
$\cala_0$ is isomorphic to the crossed product of $\cala^I_0$ with ${\mathfrak
T}^\ast$ for $\tau$. Furthermore $\cala^I_0$ is equipped with an action
$\alpha^I$ of the Lorentz group by $\ast$-automorphisms which is connected with
$\tau$ by (\ref{eq13}) and the action of the Poincar\'e group on $\cala_0$ is
given by (\ref{eq12}) and (\ref{eq15}).
\end{proposition}
In this proposition, the crossed product is taken in the category of
unital\linebreak[4] $\ast$-algebras. However one can accomodate various
functional analytic generalizations by working in the appropriate categories of
topological $\ast$-algebras. This partly depends what one has in mind for the
algebra $\calc(M)$ of functions on space-time (the linear decomposition of the
elements of $\calc(M)$ over the $u(k)$ being the Fourier transformation).\\

Let us now consider the description of the cocycle $c$ in this framework. The
restriction of $c$ to $\calc(M) \times \calc(M)$ is an $\cala_0$-valued
normalized Hochschild 2-cocycle on $\calc(M)$ which is described by the
$c(u(k),u(\ell))$, $k,\ell\in {\mathfrak T}^\ast$. The invariance and the
reality conditions yield
\begin{equation}
\alpha_{(\Lambda,a)} c(u(k),u(\ell))=c(u(\Lambda
k),u(\Lambda\ell))\exp(-i(\Lambda(k+\ell))_\mu a^\mu)
\label{eq16}
\end{equation}
\begin{equation}
c(u(k),u(\ell))^\ast=c(u(-\ell), u(-k))
\label{eq17}
\end{equation}
Define $\gamma(k,\ell)\in \cala_0$ for $k,\ell\in {\mathfrak T}^\ast$ by
\begin{equation}
c(u(k),u(\ell))=\gamma(k,\ell)u(k+\ell)
\label{eq18}
\end{equation}
It follows from (\ref{eq16})  (with $\Lambda=I$) that the $\gamma(k,\ell)$ are
translationally invariant, that is  $\gamma(k,\ell)\in \cala^I_0, \ \forall
k,\ell\in {\mathfrak T}^\ast$. The reality condition (\ref{eq17}) becomes
\begin{equation}
\gamma(k,\ell)^\ast=\tau_{k+\ell} \gamma(-\ell,-k)
\label{eq19}
\end{equation}
and (\ref{eq16}) yields the action
\begin{equation}
\alpha^I_\Lambda\gamma(k,\ell)=\gamma(\Lambda k,\Lambda \ell)
\label{eq20}
\end{equation}
On the other hand the cocycle relation on $\calc(M)$
\[
\begin{array}{lll}
u(k)c(u(\ell),u(m))& - & c(u(k)u(\ell),u(m))\\
&+ & c(u(k),u(\ell)u(m)) - c(u(k),u(\ell))u(m)=0
\end{array}
\]
is equivalent to the relation
\begin{equation}
\tau_k \gamma(\ell,m)-\gamma(k+\ell,m)+\gamma(k,\ell+m)-\gamma(k,\ell)=0
\label{eq21}
\end{equation}
and, since $\gamma(k,0)=\gamma(0,k)=0$ follows from
$c(\bbbone,x)=c(x,\bbbone)=0$, one has the following result:

\begin{lemma}
The mapping $\gamma:{\mathfrak T}^\ast\times {\mathfrak T}^\ast \rightarrow
\cala^I_0$ defined by (\ref{eq18}) is an $\cala^I_0$-valued normalized
2-cocycle of the group ${\mathfrak T}^\ast$ for $\tau$ which is real, in the
sense of (\ref{eq19}), and Lorentz-invariant, in the sense of (\ref{eq20}).
\end{lemma}

Thus starting from $(\cala^I_0, \tau, \alpha^I)$, one reconstructs $\cala_0$
with its Poincar\'e automorphism by using Proposition 2 and one reconstructs
the restriction to $\calc(M)\times \calc(M)$ of the Hochschild cocycle $c$ from
the $\cala^I_0$-valued group cocycle $\gamma$ as in Lemma 1 by formula
(\ref{eq18}). The missing items are the values of $c$ on $\cala^I_0\times
\calc(M)$ and on $\cala^I_0 \times \cala^I_0$.
Indeed let $f$ and $g$ be arbitrary elements of $\cala^I_0$ and $k,\ell$ be in
${\mathfrak T}^\ast$. Then the cocycle relation implies
\begin{equation}
\begin{array}{ll}
c(fu(k),u(\ell)g)&
=(f\gamma(k,\ell)\tau_{k+\ell}(g)+c(f,\tau_{k+\ell}(g)))u(k+\ell)\\
& + c(f\tau_{k+\ell}(g),u(k+\ell))-fu(k)c(u(\ell),g)-c(f,u(k))u(\ell)g\\
& + fc(u(k+\ell),g)-f c(\tau_{k+\ell}(g),u(k+\ell))
\end{array}
\label{eq22}
\end{equation}
This relation expresses the cocycle $c$ on $\cala_0\times \cala_0$
in terms of its restrictions to $\calc(M)\times \calc(M)$, to $\cala^I_0\times
\cala^I_0$, to $\cala^I_0\times \calc(M)$ and to $\calc(M)\times \cala^I_0$. By
taking into account the reality condition (\ref{eq17}), the restriction of $c$
to $\calc(M)\times \cala^I_0$ can be expressed in terms of its restriction to
$\cala^I_0\times \calc(M)$. Thus $c$ is known whenever its restrictions to
$\calc(M)\times \calc(M)$, to $\cala^I_0\times
\cala^I_0$ and to $\cala^I_0\times \calc(M)$ are known. It follows from
(\ref{eq6}) that the restriction $c^I$ of $c$ to
$\cala^I_0\times \cala^I_0$ is $\cala^I_0$-valued and therefore, $c^I$ is a
normalized ($\cala^I_0$-valued) 2-cocycle of $\cala^I_0$ which satisfies the
invariance condition
\begin{equation}
\alpha^I_\Lambda c^I(f,g)=c^I(\alpha^I_\Lambda f,\alpha^I_\Lambda g),
\label{eq23}
\end{equation}
and the reality condition
\begin{equation}
c^I(f,g)^\ast=c^I(g^\ast,f^\ast),
\label{eq24}
\end{equation}
for any $f,g\in \cala^I_0$ and $\Lambda\in {\mathfrak L}$.\\
We define $\lambda(f,k)$ for $f\in\cala^I_0$ and $k\in {\mathfrak T}^\ast$ by
\begin{equation}
c(f,u(k))-c(u(k),\tau_{-k}(f))=\lambda(f,k)u(k)
\label{eq25}
\end{equation}
It follows from (\ref{eq12}) that $\lambda(f,k)$ belongs to $\cala^I_0$ and
therefore, $\lambda$ is a map of $\cala^I_0\times {\mathfrak T}^\ast$ into
$\cala^I_0$, linear in the first factor, which satisfies the Lorentz-invariance
condition
\begin{equation}
\alpha^I_\Lambda \lambda(f,k)=\lambda(\alpha^I_\Lambda f,\Lambda k)
\label{eq26}
\end{equation}
the normalization conditions
\begin{equation}
\lambda(f,0)=0,\ \lambda(\bbbone,k)=0
\label{eq27}
\end{equation}
and the reality condition
\begin{equation}
\lambda(f,k)^\ast=-\tau_k(\lambda(\tau_{-k}(f^\ast),-k))
\label{eq28}
\end{equation}
The cocycle relation for $c$ implies that one has
\begin{equation}
f\gamma(k,\ell)-\gamma(k,\ell)f=\tau_k(\lambda(\tau_{-k}(f),\ell))-\lambda(f,k+\ell)+\lambda(f,k)
\label{eq29}
\end{equation}
and
\begin{equation}
c^I(f,g)-\tau_k(c^I(\tau_{-k}(f),\tau_{-k}(g)))=f\lambda(g,k)-\lambda(fg,k)+\lambda(f,k)g
\label{eq30}
\end{equation}
for any $k,\ell\in {\mathfrak T}^\ast$ and $f,g\in {\cala}^I_0$. Conversely if
$c^I,\gamma,\lambda$ as above are such that $c^I$ is a Hochschild 2-cocycle on
$\cala^I_0$ and such that (\ref{eq21}), (\ref{eq29}) and (\ref{eq30}) are
satisfied then $c$ is a 2-cocycle on $\cala_0$, i.e. the cocycle relation for
$c$ is equivalent to the cocycle relation for $c^I$ and the relations
(\ref{eq21}), (\ref{eq29}) and (\ref{eq30}).\\
We now describe another way to present these relations. We first extend the
action $\tau$ of ${\mathfrak T}^\ast$ on $\cala^I_0$ into an action on the
$\cala^I_0$-valued Hochschild cochains on $\cala^I_0$ by setting
$\tau_k(\omega)(f_1,\dots,f_n)=\tau_k(\omega(\tau_{-k}(f_1),\dots,\tau_{-k}(f_n)))$ for $k\in {\mathfrak T}^\ast$, $f_i\in \cala^I_0$, where $\omega\in C^n(\cala^I_0,\cala^I_0)$ is an $\cala^I_0$-valued Hochschild $n$-cochain on $\cala^I_0$. Let $C^{r,s}$ denote the space of $C^r(\cala^I_0,\cala^I_0)$-valued $s$-cochains on the group ${\mathfrak T}^\ast$. One has two differentials on the direct sum $C=\oplus C^{r,s}$. The first one is the composition with the Hochschild differential of $C(\cala^I_0,\cala^I_0)$ which will be denoted by $\delta^{(1,0)}$. The second one is the group differential of ${\mathfrak T}^\ast$ for its action $\tau$ on $C^r(\cala^I_0,\cala^I_0)$ which will be denoted by $(-1)^r\delta^{(0,1)}$. One has $\delta^{(1,0)}(C^{r,s})\subset C^{r+1,s}$, $\delta^{(0,1)}(C^{r,s})\subset C^{r,s+1}$ and $\delta^{(1,0)}\delta^{(0,1)}+\delta^{(0,1)}\delta^{(1,0)}=0$. Therefore $\delta=\delta^{(1,0)}+\delta^{(0,1)}$ is again a differential of degree one for the total d!
egree $r+s$. One can restrict all 
differentials to the space of normalized cochains (i.e. normalized group cochains with values in normalized Hochschild cochains). One has $c^I\in C^{2,0}$, $\lambda\in C^{1,1}$ and $\gamma \in C^{0,2}$. Furthermore, the cocycle condition for $c^I$ and the relations (\ref{eq21}), (\ref{eq29}), (\ref{eq30}) reduce to $(\delta^{(1,0)}+\delta^{(0,1)})(c^I+\lambda+\gamma)=0$. Indeed the (3,0) part of this relation $\delta^{(1,0)}c^I=0$ is the cocycle relation for $c^I$, the (2,1) part $\delta^{(0,1)}c^I+\delta^{(1,0)}\lambda=0$ is relation (\ref{eq30}), the (1,2) part $\delta^{(0,1)}\lambda+\delta^{(1,0)}\gamma=0$ is relation (\ref{eq29}) and the (0,3) part $\delta^{(0,1)}\gamma=0$ is the group cocycle relation (\ref{eq21}).

\begin{proposition}
The cocycle condition for $c$ is equivalent to
$(\delta^{(1,0)}+\delta^{(0,1)})(c^I+\lambda+\gamma)=0$, i.e. to the cocycle
condition for $c^I+\lambda+\gamma$ in $C$. Moreover the addition of the
Hochschild coboundary of a translationnally invariant 1-cochain to $c$ is
equivalent to the addition to $c^I+\lambda+\gamma$ of a term
$(\delta^{(1,0)}+\delta^{(0,1)})(\beta^{(1,0)}+\beta^{(0,1)})$ with
$\beta^{(1,0)}\in C^{(1,0)}$ and $\beta^{(0,1)}\in C^{(0,1)}$.
\end{proposition}
Indeed, if one adds to $c$ the coboundary of $b$ with $b$ such that
$\alpha_{(1,a)}(b(f))=b(\alpha_{(1,a)}(f))$ for $a\in {\mathfrak T}$ and $f\in
\cala_0$, then $\beta^{(1,0)}$ is the restriction of $b$ to $\cala^I_0$ and
$\beta^{(0,1)}$ is given by $\beta^{(0,1)}(k)=b(u(k))u(-k)$ for $k\in
{\mathfrak T}^\ast$; the invariance of $b$ implies that $\beta^{(1,0)}$ and
$\beta^{(0,1)}$ are $\cala^I_0$-valued.\\

The above result is a particular case of results of \cite{Nis} on the
Hochschild cohomology of crossed products. Here the simplification comes from
the fact that we are only interested on cochains of $\cala_0$ which are
invariant by translations (and in fact by Poincar\'e transformations).\\

It is worth noticing here that the data $c^I, \lambda$ and $\gamma$ do not
determine $c$ completely. However, by using the formula (\ref{eq22}) one can
show that $c$ is determined by $c^I,\lambda$ and $\gamma$ to within the
coboundary of a translationally invariant 1-cochain $b$ of $\cala_0$ which is
such that the corresponding $\beta^{(1,0)}$ and $\beta^{(0,1)}$ as above vanish
identically.\\

For $\kappa\not= 0$, the algebra $\cala_\kappa$ is generated by the
noncommutative version of the functions on space-time. This implies that
$\cala_0$ is generated by the algebra $\calc(M)$ of functions on space-time and
the iterated applications of the cocycle $c$ on $\calc(M)$. More precisely, we
define an increasing filtration $F^n$ of $\cala_0$ by unital $\ast$-subalgebras
$F^n(\cala_0)$ by setting $F^0(\cala_0)=\calc(M)$ and $F^{n+1}(\cala_0)=\{$ the
subalgebra of $\cala_0$ generated by $F^n(\cala_0)$ and by
$c(F^r(\cala_0),F^s(\cala_0)\}$ for $r+s=n (\in \Bbb N)$. Then, our assumption
means that one has $\cala_0=\cup_nF^n(\cala_0)=F^\infty(\cala_0)$.
Correspondingly one has an increasing filtration of $\cala^I_0$ by unital
$\ast$-subalgebras $F^n(\cala_0^I)=F^n(\cala_0)\cap \cala^I_0$ which is such
that $F^0(\cala^I_0)=\Bbb C\bbbone$, $\cala^I_0=\cup
F^n(\cala^I_0)=F^\infty(\cala^I_0)$ and which is Lorentz-invariant, i.e.
$\alpha^I_\Lambda F^n(\cala^I_0) \subset F^n(\cala^I_0)$ for $\Lambda \in
\mathfrak L$.

\section{Discussion}

The simplest cases correspond to the situation where $\calc(M)$ is in the
center of $\cala_0$, i.e. where $\tau$ is trivial. In this case, one can assume
that, up to a coboundary, $c$ is antisymmetric on $\calc(M)$ and in fact on
$Z(\cala_0$), (this corresponds to a weak regularity assumption on the
commutative algebra $Z(\cala_0)$). It then follows from the discussion in the
end of last Section that one has $\cala_0=Z(\cala_0)$ i.e. that $\cala_0$ is a
commutative Poisson algebra and that the $\cala_\kappa$ are obtained by
``quantization" of $\cala_0$. In this case $\cala_0$ is, as algebra, the tensor
product $C(M)\otimes \cala^I_0$. Since the $\{x^\mu,x^\nu\}$ are elements of
$\cala^I_0$ and since the Lorentz group acts there by automorphisms,
$\cala^I_0$ must contain as subalgebra an algebra of functions on a union of
Lorentz orbits of antisymmetric 2-tensors. The $\{x^\mu,x^\nu\}$ are then
coordinates functions on this manifold of antisymmetric tensors and the orbits
occurring are labeled by the pairs $(\alpha,\beta)$ of real numbers such that
$\alpha$ is in the spectrum of
$g_{\mu\mu'}g_{\nu\nu'}\{x^\mu,x^\nu\}\{x^{\mu'},x^{\nu'}\}$ and $\beta$ is in
the spectrum of
$\varepsilon_{\alpha\beta\gamma\delta}\{x^\alpha,x^\beta\}\{x^\gamma,x^\delta\}$. It is natural to require time reversal and  parity be defined and therefore to assume that whenever one has the orbit $(\alpha,\beta)$, one also has the orbit $(\alpha,-\beta)$. When furthermore one has $\{x^\lambda,\{x^\mu,x^\nu\}\}=0$\ $(\forall \lambda,\mu,\nu)$, then $\cala^I_0$ is just such an algebra of functions on a union of orbits of antisymmetric tensors. This is precisely the case for the algebra $\cala_0$ which is the $\kappa=0$ limit of the model of Doplicher, Fredenhagen and Robert \cite{DFR} where the orbits there are (0,1) and (0,-1). It is not very difficult to construct examples with $\{x^\lambda,\{x^\mu,x^\nu\}\}\not=0$.\\

As pointed out above, in order to have a noncommutative algebra $\cala_0$, it
is necessary (and obviously sufficient) that the algebra $\calc(M)$ of
functions on space-time be not included in the center of $\cala_0$, which means
that $\tau$ is non trivial. Since by assumption the cocycle $\gamma$ defined by
(\ref{eq18}) is non trivial, the simplest cases with $\tau$ non trivial are the
cases where $\lambda$ and $c^I$ vanish. In such a case, it follows from
(\ref{eq29}) that the image of $\gamma$ is in the center of $\cala^I_0$ and
that $F^{n+1}(\cala_0)=F^1(\cala_0)$, $\forall n\in \mathbb N$. Therefore, in
view of the discussion in the end of last section, $\cala_0=F^1(\cala_0)$ and
$\cala^I_0$ is the commutative algebra generated by the $\gamma(k,\ell)$, (or,
equivalently, by the  $\{x^\mu,x^\nu\}$). Thus in such a case $\cala_0$ is the
crossed product of the commutative algebra $\cala^I_0$ with ${\mathfrak
T}^\ast$ for $\tau$. An example of this situation is provided by the $\kappa=0$
limit of an example elaborated by Doplicher and Fredenhagen described in
Section 2 of \cite{Dop}. It is worth noticing here that in the $\kappa=0$
limits of examples of \cite{DFR} and \cite{Dop} the orbits of antisymmetric
2-tensors occuring are $(0,1)$ and $(0,-1)$, (recall that in these examples
$\cala^I_0$ is commutative). This is connected with the fact that these authors
construct $\cala_\kappa$ in such a way that physically-motivated spacetime
uncertainty relations are implemented.\\

Generically $\tau$ is non-trivial and $\lambda$ and $c^I$ do not vanish. A
simple example of this kind can be easily found. For $\kappa\in \mathbb R$, let
$\cala_\kappa$ be the unital $\ast$-algebra generated by hermitian elements
$x^\mu, L^\mu, I^{\mu\nu}\ (\mu,\nu=0,1,2,3)$ satisfying the relations
\begin{eqnarray}
[x^\mu,x^\nu]& = &i\kappa I^{\mu\nu}\nonumber\\
{[x^\lambda,I^{\mu\nu}]}& = &i(g^{\lambda\nu}L^\mu-g^{\lambda\mu}L^\nu)\\
{[x^\mu,L^\nu]}& = &i\kappa I^{\mu\nu}\nonumber
\end{eqnarray}

\begin{eqnarray}
[I^{\lambda\rho},I^{\mu\nu}]& = &
i(g^{\lambda\nu}I^{\mu\rho}-g^{\rho\nu}I^{\mu\lambda}+g^{\rho\mu}I^{\nu\lambda}-g^{\lambda\mu}I^{\nu\rho})\nonumber\\
{[L^\lambda,I^{\mu\nu}]}& = & i(g^{\lambda\nu}L^\mu-g^{\lambda\mu}L^\nu)\\
{[L^\mu,L^\nu]}& = & i\kappa I^{\mu\nu}\nonumber
\end{eqnarray}
where $g^{\mu\nu}$ denotes the minkowskian metric with
$g^{00}=-1,g^{11}=g^{22}=g^{33}=1$. It follows from these relations that for
$\kappa\not=0$, $\cala_\kappa$ is generated by the $x^\mu$. This implies that
$\cala_0$ satisfies the property
$\cala_0=\cup_nF^n(\cala_0)$ of the end of last section. In fact here one has
again $F^{n+1}(\cala_0)=F^1(\cala_0)$, $\forall n\in \mathbb N$. Furthermore
there is (for any $\kappa$) an action of the Poincar\'e group $\mathfrak P$ by
$\ast$-automorphisms $(\Lambda,a)\mapsto \alpha_{(\Lambda,a)}$ on
$\cala_\kappa$ given by
\[
\alpha_{(\Lambda,a)}x^\mu=\Lambda^{-1\mu}_\nu(x^\nu-a^\nu\bbbone),\ \
\alpha_{(\Lambda,a)}L^\mu=\Lambda^{-1\mu}_\nu L^\nu,\ \
\alpha_{(\Lambda,a)}I^{\mu\nu}=\Lambda^{-1\mu}_\alpha\Lambda^{-1\nu}_\beta
I^{\alpha\beta}
\]
The commutation relations (32) between the $I^{\mu\nu}$ and the $L^\lambda$ are
the relations of the Lie algebra of $SO(4,1)$ if $\kappa>0$, of $SO(3,2)$ if
$\kappa<0$ and of the Poincar\'e group if $\kappa=0$. It follows that the
$I^{\mu\nu}$ and the $L^\lambda$ generate the corresponding enveloping algebra.
The $x^\mu-L^\mu$ are in the center $Z(\cala_\kappa)$ of $\cala_\kappa$.
Therefore the algebra $\cala_\kappa$ is the tensor product of the commutative
algebra generated by the $x^\mu-L^\mu\ (\mu=0,1,2,3)$ with the above enveloping
algebra. In fact the center $Z(\cala_\kappa)$ of $\cala_\kappa$ is generated by
the $x^\mu-L^\mu$ and the two casimirs  $\mathfrak c_2$ and $\mathfrak c_4$
given by
\begin{equation}
{\mathfrak c}_2=\kappa g_{\lambda\mu} g_{\rho\nu} I^{\lambda \rho}I^{\mu\nu}+2
g_{\alpha\beta} L^\alpha L^\beta
\end{equation}
\begin{equation}
{\mathfrak c}_4 = g^{\rho\rho'}\varepsilon_{\rho\lambda\mu\nu}L^\lambda
I^{\mu\nu}\varepsilon_{\rho'\lambda'\mu'\nu'}L^{\lambda'}I^{\mu'\nu'}
\end{equation}
where $\varepsilon_{\rho\lambda\mu\nu}$ is the completely antisymmetric tensor
with $\varepsilon_{0123}=1$. Therefore $\cala_0$ is the tensor product of the
commutative algebra generated by the $x^\mu-L^\mu$ with the enveloping algebra
of the Poincar\'e Lie algebra generated by the $I^{\mu\nu}$ and the
$L^\lambda$. Notice however that the $I^{\mu\nu}$ and the $L^\lambda$ have
nothing to do with the physical Poincar\'e group $\mathfrak P$ acting on
$\cala_0$ by the automorphisms $\alpha_{(\Lambda,a)}$; in particular the
$L^\lambda$ have the dimension of a length. The algebra $\cala^I_0$ is simply
the subalgebra generated by the $I^{\mu\nu}$ and the $L^\lambda$ ($\simeq$
enveloping algebra of a Poincar\'e Lie algebra). The $\tau_k$ are given by
$\tau_k(I^{\mu\nu})=I^{\mu\nu}+k^\mu L^\nu-k^\nu L^\mu$ and
$\tau_k(L^\lambda)=L^\lambda$. One can compute the cocycle $c$, e.g. one has
$c(x^\mu,x^\nu)=-\frac{1}{2} I^{\mu\nu}$. Finally, allowing exponentials, the
cocycle $\gamma$ is given by
\[
\gamma(k,\ell)=-\frac{i}{2}\left(k_\mu\ell_\nu
I^{\mu\nu}-\left(\frac{2}{3}k_\rho\ell^\rho+\frac{1}{3}\ell_\rho\ell^\rho\right)k_\lambda L^\lambda+\left(\frac{1}{3}k_\rho\ell^\rho+\frac{2}{3}k_\rho k^\rho\right)\ell_\lambda L^\lambda\right)
\]
and one has $\lambda(L^\mu,k)=-k_\nu\left(I^{\mu\nu}+\frac{1}{2}(k^\mu
L^\nu-k^\nu L^\mu)\right)$ and $c^I(L^\mu,L^\nu)=-\frac{1}{2}I^{\mu\nu}$ and so
forth. Since $\mathfrak c_2$ and $\mathfrak c_4$ are in the center of
$\cala_\kappa$ and since they are translationally invariant it is natural to
fix them and to add the corresponding relations in the definition of the
$\cala_\kappa$. Since the element $\mathfrak c_2$ has the dimension of a length
squared, there are three natural ways to fix it: (i) $\mathfrak c_2=\kappa$,
(ii) $\mathfrak c_2=-\kappa$, (iii) $\mathfrak c_2=0$. Any of these choices
leads to $g_{\mu\nu}L^\mu L^\nu=0$ in $\cala_0$ (in fact in $\cala^I_0$).
Remembering that $\cala^I_0$ has the structure of the enveloping algebra of a
Poincar\'e Lie algebra, the latter condition is the analogue of a zero mass
condition (although here it has a different meaning as we pointed out). Thus if
one assumes furthermore that $\mathfrak c_4$ is fixed in such a way that the
representations occuring are ``zero mass" representations with strictly
positive ``spin", one finds again a characteristic two sheeted structure for
$\cala^I_0$ (by allowing the two helicities). The origin of the frequent
occurence of this two sheeted structure in this context is obviously due to the
fact that the full Lorentz group is not connected.\\

In all the above examples the $\cala_\kappa$ are Poincar\'e-covariant. This is
not needed, only $\cala_0$ must have this property. In fact one can easily
produce an example where only $\cala_0$ is Poincar\'e covariant by deforming
the universal enveloping algebras occuring in the previous example for
$\kappa\not= 0$ (in the sense used in the theory of quantum groups).\\

If one wishes to establish a connection between the extra factor $\cala^I_0$
here and the one occuring in recent noncommutative versions of gauge theory
(\cite{dvkm}, \cite{C}, \cite{coq}), one should remember that, according to our
analysis, $\cala^I_0$ must be infinite dimensional and that it can be
noncommutative only if some of its elements do not commute with the functions
on space-time. Concerning the first point, one could expect that, by some
contractibility argument, a finite dimensional approximation can be found.
However the second point remains. This suggests that it is worth trying to
enlarge the setting of the noncommutative models of gauge theory.

\section*{Acknowledgements}

It is a pleasure to thank Ivan Todorov for stimulating discussions, Victor
Nistor for enlightening correspondence and Sunggoo Cho for pointing out to us a mistake in the previous version of the paper.


\newpage

\begin{thebibliography}{99}


\bibitem[1]{BFLS} F. Bayen, M. Flato, C. Fronsdal, A. Lichnerowicz, D.
Sternheimer. \newblock Deformation theory and quantization. I. Deformations of
symplectic structures.
\newblock {\em Ann. Phys.\/ \bf 111} (1978) 61 and
\newblock II. Physical applications.
\newblock {\em Ann. Phys.\/ \bf 111} (1978) 111.

\bibitem[2] {B} C. Bertrand.
\newblock Calcul diff\'erentiel non commutatif sur une alg\`ebre de Hopf. Un
outil essentiel, la repr\'esentation $\tau$.
\newblock {\em Bull. Sci. math.\/ \bf 119} (1995) 49.

\bibitem[3] {C} A. Connes, J. Lott.
\newblock Particles models and non commutative geometry.
\newblock {\sl Nucl. Phys. B (Proc. Suppl.)\bf 11B} (1990) 19.\\
\newblock A. Connes.
\newblock Non-commutative geometry.
\newblock Academic Press, 1994.

\bibitem[4]{coq}
\newblock R. Coquereaux.
\newblock
Higgs fields and superconnections,
in Differential Geometric Methods in Theoretical Physics, Rapallo 1990 (C.
Bartocci, U. Bruzzo, R. Cianci, eds),
{\em Lecture Notes in Physics\/ \bf 375}, Springer Verlag 1991.\\
\newblock R. Coquereaux, R. H\"au\ss ling, F. Scheck.
Algebraic connections on parallel universes,
\newblock {\em Int. J. Mod. Phys.\/ \bf A10} (1995) 89-98.

\bibitem[5]{DKP}
\newblock C. De Concini, V.G. Kac, C. Procesi,
\newblock Some quantum analogues of solvable Lie groups.
\newblock {\sl Geometry and Analysis}, Tata Institute of Fundamental Research,
hep-th/9308138.


\bibitem [6]{DS}
\newblock S. Deser
\newblock General relativity and the divergence problem in quantum field
theory.
\newblock {\em Rev. Mod. Phys.\/ \bf29} (1957) 417.

\bibitem [7] {Dop}
\newblock S. Doplicher.
\newblock Quantum spacetime.
\newblock {\em Ann. Inst. Henri Poincar\'e\/ \bf 64} (1996) 543.

\bibitem [8]{DFR}
\newblock S. Doplicher, K. Fredenhagen, J.E. Roberts.
\newblock The quantum structure of spacetime at the Planck scale and quantum
fields.
\newblock {\em Commun. Math. Phys.\/ \bf 172} (1995) 187.

\bibitem[9]{dvkm}
\newblock M. Dubois-Violette, R. Kerner, J. Madore.
\newblock
Gauge bosons in a non-commutative geometry,
\newblock
{\em Phys. Lett.\/ \bf B217} (1989) 485.\\
\newblock
M. Dubois-Violette, R. Kerner, J. Madore.
\newblock
Classical bosons in a non-commutative geometry,
\newblock {\em Class. Quantum Grav.\/ \bf 6} (1989) 1709.\\
\newblock
M. Dubois-Violette, R. Kerner, J. Madore
\newblock
Non-commutative differential geometry of matrix algebras,
\newblock {\em J. Math. Phys.\/ \bf 31} (1990) 316.\\
\newblock
M. Dubois-Violette, R. Kerner, J. Madore.
\newblock
Non-commutative differential geometry and new models of gauge theory,
\newblock
{\em J. Math. Phys.\/ \bf 31} (1990) 323.

\bibitem[10]{G}
\newblock M. Gerstenhaber.
\newblock On the deformation of rings and algebras.
\newblock {\em Ann. Math.\/ \bf 79} (1964) 59.

\bibitem[11] {Nis}
\newblock V. Nistor.
\newblock Group cohomology and the cyclic cohomology of crossed products.
\newblock{\em Invent. Math.\/ \bf 99} (1990) 411.

\bibitem[12]{RVW}
\newblock N. Reshetikhin, A. Voronov, A. Weinstein.
\newblock Semiquantum geometry.
\newblock Preprint Juin 1996.

\bibitem[13]{W}
S.L.Woronowicz,
Differential calculus on compact matrix pseudogroups (quantum groups),
\newblock
{\em Commun. Math. Phys.\/ \bf 122} (1989) 125.


\end{thebibliography}
\end{document}